\documentclass[11pt,a4paper,reqno]{amsart}
\usepackage{amsmath,amssymb,amsfonts,amsthm,mathtools}
\usepackage{enumerate}
\usepackage[colorlinks=true,linkcolor=black,citecolor=black,urlcolor=black,pagebackref=false]{hyperref}

\numberwithin{equation}{section}
\allowdisplaybreaks

\usepackage{ulem}
\newcommand{\R}{\mathbb R}

\newcommand{\supp}{\mathrm{supp}}
\usepackage{booktabs}
\newcommand{\one}{\mathbf 1}

\newcommand{\Z}{\mathbb{Z}}
\newcommand{\T}{\mathbb{T}}
\newcommand{\1}{\mathbf{1}}
\newcommand{\dd}{\mathrm{d}}
\newcommand{\FT}[1]{\widehat{#1}}
\makeatletter
\@namedef{subjclassname@2020}{\textup{2020} Mathematics Subject Classification}
\makeatother
\newtheorem{thm}{Theorem}[section]
\newtheorem{prop}[thm]{Proposition}
\newtheorem{lem}[thm]{Lemma}
\newtheorem{coro}[thm]{Corollary}

\newtheorem{rem}[thm]{Remark}
\newtheorem{ques}[thm]{Question}
\newtheorem{conj}[thm]{Conjecture}

\title[On the spectrality]
{On the spectrality of non-homogeneous self-similar measures}

\author{Yi-Qiu Mao}
\address{Y.-Q. Mao, School of Mathematics and Information Science, Guangzhou University, Guangzhou, 510006, P.~R.~China}
\email{yqmao@gzhu.edu.cn}

\author{Zhi-Yi Wu}
\address[Zhi-Yi Wu]{School of Mathematics and Information Science, Guangzhou University, Guangzhou, 510006, P.~R.~China}
\email{zhiyi\_wu2021@163.com}


\begin{document}
\begin{abstract}
We develop a finite-cylinder difference-set obstruction to the existence of
an exponential orthonormal basis for a compactly supported probability
measure.  This obstruction uses the equal-scale and equal-mass structure
of finitely many disjoint affine copies of the measure, together with
Parseval's identity and an elementary atom bound for positive measures.
We apply it to prove that the golden-mean self-similar measure generated
by $x/2$ and $(x+2)/4$ is non-spectral.  Together with a weight-rigidity
argument, this excludes spectrality for every positive choice of weights
in this iterated function system.  We also analyze a logarithmically incommensurable
example generated by $x/2$ and $(x+2)/3$, proving non-spectrality for
all first-branch weights $p\in[1/3,5/6]$, including the natural weights.
The latter proof reduces to two finite, exactly checkable integer
difference-set conditions; an independent verification program is included in the paper.
\end{abstract}

\maketitle

\section{Introduction}
A Borel probability measure \(\mu\) on \(\mathbb{R}^d\) is called a {\it spectral measure} if there exists a countable set \(\Lambda \subseteq \mathbb{R}^d\) such that \(E(\Lambda):=\{e_\lambda(x):\lambda\in \Lambda\}\) with $e_\lambda(x)=e^{2\pi i\lambda \cdot x}$, is an orthonormal basis for \(L^2(\mu)\), and \(\Lambda\) is called a {\it spectrum} of \(\mu\). The study of spectral measures can be traced back to the celebrated Fuglede conjecture in 1974, which asserted that the normalised Lebesgue measure restricted on a bounded measurable set \(\Omega\) is spectral if and only if $\Omega$ tiles \(\mathbb{R}^d\) by translations \cite{Fug74}. Although this conjecture was later disproved in higher dimensions \cite{Tao04}, it sparked a profound interest in the interplay between orthogonal exponentials and geometric structures.

In 1998, Jorgensen and Pedersen \cite{JP98} constructed the first example of a non-atomic singular spectral measure, i.e., the standard middle-fourth Cantor measure $\mu_{4,\{0,2\}}$. They also proved in the same paper that the standard middle-third Cantor measure $\mu_{3,\{0,2\}}$ is not a spectral measure. Following this discovery, the spectrality of fractal measures including self-similar measures, self-affine measures, and Moran measures has been extensively explored; see \cite{AFL19,Dai12,DC21,DHLai19,DL14,LW02} and the references therein. All these known spectral measures admit infinite convolution of finitely supported discrete measures, which provide a powerful route to explicit spectra.  However, not every fractal measure admits such an infinite-convolution structure. In particular, it need not be available for iterated function systems with unequal contraction ratios, and our arguments below do not rely on it.

In this paper, we consider a class of non-homogeneous self-similar measures, which is studied by Saglietti, Shmerkin and  Solomyak in \cite{SSS18}. Let $r_1,r_2\in(0,1)$, $p\in(0,1)$ and $t_1,t_2\in\R$.
Consider the iterated function system (IFS): $f_1(x)=r_1x+t_1$ and $f_2(x)=r_2x+t_2$.
By Hutchinson's theorem \cite{Hut81}, there exists a unique Borel
probability measure $\mu_p$ with nonempty compact support $K$
satisfying
\begin{align}\label{def:self-s}
\mu_p(E) = p\,\mu_p\bigl(f_{1}^{-1}(E)\bigr) + (1 - p)\,\mu_p\bigl(f_{2}^{-1}(E)\bigr),
\end{align}
for any Borel set $E\subseteq \mathbb{R}$ and $K$ is the unique compact set satisfying
\begin{align}\label{defadd:self-sset}
K=f_{1}(K)\cup f_{2}(K).
\end{align}
The probability measure $\mu_p$ satisfying \eqref{def:self-s} is
called a {\it self-similar measure}, and $K$ is called a {\it self-similar} set or the attractor of the IFS $\{f_1,f_2\}$. We say that the IFS $\{f_1,f_2\}$ satisfies the {\it strong separation condition}
(SSC) if $f_1(K)\cap f_2(K)=\varnothing$.

Our original motivation was to construct a class of non-homogeneous self-similar spectral measures. As a first step in this direction, in this paper we obtain the following necessary condition for spectrality.

\begin{thm}\label{thm:mainintro0}
Let $\mu_p$ be the measure defined in \eqref{def:self-s}. Suppose $\{f_1,f_2\}$ satisfies the strong separation condition, and \(\log r_1/\log r_2\in \mathbb{Q}.\)
If $\mu_p$ is a spectral measure, then $p=r_1^s$, where $s$ is the unique positive solution of
\(r_1^s+r_2^s=1.\)
\end{thm}

Theorem~\ref{thm:mainintro0} gives a weight restriction under
logarithmic commensurability.  It is natural to ask whether these
weights suffice for spectrality.  We answer this negatively for an
explicit strongly separated, non-homogeneous system. Specifically, consider the iterated function system
\[
f_1(x) = \frac{x}{2}, \qquad f_2(x) = \frac{x+2}{4},
\]
with weight \(p = \frac{\sqrt{5}-1}{2}\). In this case, the self-similar measure is called the {\it golden-mean self-similar measure}, and we denote it by \(\mu_g\). Since $p+p^2=1$, the natural weights are $p=(1/2)^s$ and $p^2=(1/4)^s$ for $s=\log_2(1/p)>0$.  Thus they satisfy
Theorem~\ref{thm:mainintro0}, but the following theorem rules out
spectrality.

\begin{thm}\label{thm:mainintro1}
The golden-mean self-similar measure \(\mu_g\) is not a spectral measure.
\end{thm}

Combining Theorems \ref{thm:mainintro0} and \ref{thm:mainintro1}, we immediately have the following corollary.
\begin{coro}\label{cor:golden-all-weights}
For the IFS $f_1(x)=x/2$, $f_2(x)=(x+2)/4$, the self-similar
measure $\mu_p$ is non-spectral for every $p\in(0,1)$.
\end{coro}

For the proof of Theorem \ref{thm:mainintro1}, our method is a finite-scale obstruction built from disjoint
cylinders with a common contraction and mass.  Under a spectrality
hypothesis, Parseval's identity turns the arithmetic of their
translation differences into vanishing Fourier moments of a positive
measure on the circle.  A Fej\'{e}r-kernel bound then conflicts with
an atom whose mass is supplied by the zero frequency.  This approach
neither assumes a convolution representation nor presupposes any zero-free property of the original Fourier transform, in contrast to the lattice-containment assumption on the Fourier zero set used in \cite{AW21}.

The cylinder obstruction is not confined to logarithmically
commensurable contraction ratios.  To illustrate this point, let
\[
  g_1(x)=\frac{x}{2},\qquad g_2(x)=\frac{x+2}{3},
\]
and denote by $\nu_p$ their self-similar measure with weights
$p$ and $1-p$.  Here $\log(1/2)/\log(1/3)\notin\mathbb Q$.
Equal-scale cylinders nevertheless arise from words with the same
numbers of the two symbols.  Their translation differences lead to
the following interval of non-spectral parameters.

\begin{thm}\label{thm:incommensurable-interval}
The self-similar measure $\nu_p$ generated by $g_1(x)=x/2$ and
$g_2(x)=(x+2)/3$ is not spectral for every
\(
  p\in\left[\frac13,\frac56\right].
\)
\end{thm}

\begin{coro}\label{cor:incommensurable-natural}
Let $s>0$ be determined by $2^{-s}+3^{-s}=1$.  The self-similar
measure generated by $g_1,g_2$ with the natural weights
$p=2^{-s}$ and $1-p=3^{-s}$ is not spectral.
\end{coro}

The finite integer-difference verifications used for
Theorem~\ref{thm:incommensurable-interval} are specified by an
exact recurrence, an independently reproducible verification
algorithm, and complete witness data supplied with this manuscript.
No floating-point computation of a Fourier transform enters the proof.

This paper is organized as follows. Section 2 develops the general cylinder criterion and proves
Theorem~\ref{thm:mainintro0}. Sections 3 and 4 study the golden-mean self-similar measure. Section 5 discusses the logarithmically incommensurable example.  Finally, we conclude in Section 6 with some remarks and open questions.
\section{A general finite-cylinder obstruction}
Write $\T=\{z\in\mathbb C:|z|=1\}$. We first obtain an elementary moment bound.

\begin{lem}\label{lem:atom}
Let $\sigma$ be a Borel probability measure on $\T$. Suppose that,
for an integer $N\ge1$,
\begin{equation}\label{eq:moments}
  \int_{\T}z^k\,\dd\sigma(z)=0\qquad (k=1,\ldots,N).
\end{equation}
Then
\begin{equation}\label{eq:atom_bound}
  \sigma(\{1\})\le\frac1{N+1}.
\end{equation}
\end{lem}

\begin{proof}
Since $\sigma$ is a positive measure, the complex conjugates of
\eqref{eq:moments} also vanish; thus the integral of $z^k$ is
zero whenever $0<|k|\le N$. Define
\(
  P_N(z)=\sum_{j=0}^{N}z^j.
\)
Then
\begin{align}\label{eq:geqee0}
 \int_{\T}|P_N(z)|^2\,\dd\sigma(z)
 &=\sum_{j=0}^{N}\sum_{\ell=0}^{N}
     \int_{\T}z^{j-\ell}\,\dd\sigma(z) \nonumber\\
 &=\sum_{j=0}^{N}\int_{\T}1\,\dd\sigma(z)
   =N+1.
\end{align}
On the other hand $P_N(1)=N+1$, and the integrand \(|P_N(z)|^2\) is nonnegative. Consequently
\begin{align}\label{eq:geqee1}
\int_{\T}|P_N(z)|^2\,\dd\sigma(z) \ge \int_{\{1\}}|P_N(z)|^2\,\dd\sigma(z) = |P_N(1)|^2\sigma(\{1\}).
\end{align}

Combining \eqref{eq:geqee0} and \eqref{eq:geqee1}, we obtain \eqref{eq:atom_bound}.
\end{proof}

For a finite Borel measure \(\nu\) on \(\R^d\), its Fourier transform
is defined by
\[
  \FT\nu(\xi)
  :=\int_{\R^d} e^{-2\pi i\,\xi\cdot x}\,\dd\nu(x),
  \qquad \xi\in\R^d.
\]

\begin{prop}\label{prop:obstruction}
Let $\nu$ be a compactly supported Borel probability measure on $\R$
and set $K_\nu=\supp\nu$. Let $A\subseteq\Z$ be a finite set, and let
$r>0$, $h\ne0$, $t\in\R$, $m>0$. For each $a\in A$, define
\[
  \phi_a(x)=rx+t+ha\quad\text{and}\quad K_a:=\phi_a(K_\nu).
\]
Assume that:
\begin{enumerate}[(i)]
  \item[\textup{(1)}] the sets $K_a$ $(a\in A)$ are pairwise disjoint.
  \item[\textup{(2)}] as measures on $\R$,
  \begin{equation}\label{eq:restricted}
  \nu\big|_{K_a}(\cdot)=m\,\nu(\phi_a^{-1}(\cdot))
       \qquad(a\in A).
\end{equation}
\item[\textup{(3)}] there exists some positive integer $N$ such that
\begin{equation}\label{eq:differences_general}
  \{-N,-N+1,\ldots,N-1,N\}\subseteq A-A.
\end{equation}
\end{enumerate}
If $\nu$ is a spectral measure, then
\begin{equation}\label{eq:necessary}
  m\le\frac1{N+1}.
\end{equation}
\end{prop}
\begin{proof}
Assume that $\nu$ is a spectral measure with a spectrum $\Lambda$. Since spectra are invariant under translations, without loss of generality, we assume $0\in \Lambda$.

For $a\in A$ define
\[
  g_a=m^{-1/2}\,\1_{K_a}\in L^2(\nu),
\]
where $\1_E$ denotes the indicator function of a
set $E$. Taking total masses in \eqref{eq:restricted} shows
$\nu(K_a)=m$. Pairwise disjointness therefore implies
\begin{equation*}%\label{eq:gON}
  \langle g_a,g_b\rangle_{L^2(\nu)}=\delta_{ab},
\end{equation*}
where \(\delta_{ab}\) denotes the Kronecker delta,
that is,
\[
  \delta_{ab}=
  \begin{cases}
    1, & a=b,\\
    0, & a\ne b.
  \end{cases}
\]
We use the convention
$\langle u,v\rangle=\int u\overline v\,\dd\nu$
(linear in the first argument). From \eqref{eq:restricted},
for every $\lambda\in\Lambda$,
\begin{align}
  \langle g_a,e_\lambda\rangle
  &=\frac1{\sqrt m}\int_{K_a} e^{-2\pi i\lambda x}
                              \,\dd\nu(x)\notag\\
  &=\sqrt m\int_{K_\nu}e^{-2\pi i\lambda\phi_a(y)}
                              \,\dd\nu(y)\notag\\
  &=\sqrt m\,e^{-2\pi i\lambda(t+ha)}
                  \FT\nu(r\lambda).
  \label{eq:innerproduct}
\end{align}
Because $E(\Lambda)$ is a complete orthonormal system,
Parseval's identity applied to the pair $g_a,g_b$ gives
\begin{align}
  \delta_{ab}
  &=\sum_{\lambda\in\Lambda}
      \langle g_a,e_\lambda\rangle
      \overline{\langle g_b,e_\lambda\rangle}\notag\\
  &=\sum_{\lambda\in\Lambda}
      m\,|\FT\nu(r\lambda)|^2
      e^{-2\pi ih\lambda(a-b)}.
  \label{eq:parseval}
\end{align}
The series is well-defined: for each $a$, Parseval gives
$\sum_{\lambda\in\Lambda}|
\langle g_a,e_\lambda\rangle|^2=1$;
absolute convergence of the cross-series also follows by
Cauchy--Schwarz.

Define a positive measure on the circle by
\begin{equation}\label{eq:sigma_general}
  \sigma=\sum_{\lambda\in\Lambda}
      m\,|\FT\nu(r\lambda)|^2
      \delta_{\exp(2\pi ih\lambda)},
\end{equation}
where $\delta_z$ denotes the Dirac measure at $z$. Here repeated points on $\T$ are allowed: their nonnegative
weights are added. Setting $a=b$ in \eqref{eq:parseval} yields
\begin{equation}\label{eq:sigma_prob}
  \sigma(\T)
  =\sum_{\lambda\in\Lambda}
      m\,|\FT\nu(r\lambda)|^2=1.
\end{equation}
Thus $\sigma$ is a probability measure. If
$k=b-a$ with $a,b\in A$ and $k\ne0$, then
\eqref{eq:parseval} gives
\begin{equation}\label{eq:sigma_zero}
  \int_{\T}z^k\,\dd\sigma(z)
     =\sum_{\lambda\in\Lambda}
       m\,|\FT\nu(r\lambda)|^2e^{2\pi ih\lambda k}
     =0.
\end{equation}
Hypothesis \eqref{eq:differences_general} therefore supplies all
vanishing moments of orders $1,\ldots,N$.

Finally, $0\in \Lambda$ and $\FT\nu(0)=1$ show that
the $\lambda=0$ term in \eqref{eq:sigma_general} places mass
exactly $m$ at $z=1$. Other terms may add more mass, so
\begin{equation}\label{eq:atom_at_one}
  \sigma(\{1\})\ge m.
\end{equation}
Apply Lemma~\ref{lem:atom} to \eqref{eq:sigma_zero} and combine
with \eqref{eq:atom_at_one} to conclude
$m\le\sigma(\{1\})\le1/(N+1)$.
\end{proof}

The requirement \(A\subseteq\mathbb Z\) in
Proposition~\ref{prop:obstruction} is not essential: only the
difference set \(A-A\) matters. The following real-difference
version makes this precise.

\begin{prop}\label{prop:real-difference}
The integer-lattice assumption is not essential. Let $\nu$ be a
compactly supported Borel probability measure on $\R$, and set
$K_\nu=\supp\nu$. Let $A\subseteq\R$ be a finite set, and let
$r>0$, $t\in\R$, $m>0$. For each $a\in A$, define
\[
  \phi_a(x):=rx+t+a,\qquad K_a:=\phi_a(K_\nu).
\]
Assume that:
\begin{enumerate}
 \item[\textup{(1)}] the sets $K_a$ $(a\in A)$ are pairwise disjoint;
 \item[\textup{(2)}] as measures on $\R$,
  \begin{equation}\label{eq:real-restricted}
    \nu\big|_{K_a}(\cdot)=m\,\nu(\phi_a^{-1}(\cdot))
    \qquad(a\in A).
  \end{equation}
  \item[\textup{(3)}] there exists a nonempty finite set $T\subseteq\R$ satisfying
$T-T\subseteq A-A$.
\end{enumerate}
If $\nu$ is spectral, then
\[
  m\,\#T\le1,
\]
where \(\#\) denotes the cardinality of a finite set.
\end{prop}
\begin{proof}
Let $\Lambda$ be a spectrum translated to contain $0$.
Exactly as in \eqref{eq:innerproduct}--\eqref{eq:parseval},
\[
  \eta=\sum_{\lambda\in\Lambda}
       m|\widehat\nu(r\lambda)|^2\delta_\lambda
\]
on $\R$ is a probability measure and satisfies
\[
  \int_\R e^{2\pi i(a-b)\lambda}\,\dd\eta(\lambda)
     =\delta_{ab}\qquad (a,b\in A).
\]
In particular, its Fourier moments vanish at every nonzero
element of $A-A$.  Since $T-T\subseteq A-A$, expansion of a
finite squared sum yields
\[
  \int_\R\left|\sum_{t\in T}e^{2\pi i t\lambda}\right|^2
               \dd\eta(\lambda)=\#T.
\]
But $\eta(\{0\})\ge m$, because the term $\lambda=0$ has
weight $m|\widehat\nu(0)|^2=m$.  Consequently
\[
  \#T\ge\eta(\{0\})(\#T)^2\ge m(\#T)^2,
\]
which proves the result.
\end{proof}

By Proposition~\ref{prop:real-difference}, we immediately have the following corollary.
\begin{coro}
\label{cor:capacity}
For a nonempty finite set $A\subseteq\R$ define
\[
  \kappa(A)=\max\big\{\#T:T\subseteq\R\text{ finite},\
                         T-T\subseteq A-A\big\}.
\]
This number is finite: translating a candidate $T$ so that
$0\in T$ gives $T\subseteq A-A$.  Under the hypotheses of
Proposition~\ref{prop:real-difference}, spectrality implies
\[
  m\,\kappa(A)\le1.
\]
\end{coro}
%\begin{proof}
%Let \(T\subseteq\R\) be a finite set attaining the maximum,
%so that \(\#T=\kappa(A)\) and \(T-T\subseteq A-A\).
%By Proposition~\ref{prop:real-difference}, spectrality implies
%\(m\#T\le1\), which is exactly \(m\kappa(A)\le1\).
%\end{proof}

\begin{lem}
\label{lem:mass-rigidity}
Let $\nu$ be a compactly supported probability measure on $\R$
which is spectral.  Suppose two disjoint sets $\phi_i(K_\nu)$, $i=1,2$, satisfy
\[
  \nu\big|_{\phi_i(K_\nu)}(\cdot)=m_i\nu(\phi_i^{-1}(\cdot)),
  \qquad\phi_i(x)=rx+t_i,\quad m_i>0,
\]
with the same nonzero contraction $r$.  Then $m_1=m_2$.
\end{lem}

\begin{proof}
For each $i$, the normalized indicator
$h_i=m_i^{-1/2}\one_{\phi_i(K_\nu)}$ has $L^2(\nu)$-norm one.
For a spectrum $\Lambda$, Parseval's identity gives
\[
  1=\sum_{\lambda\in\Lambda}
       |\langle h_i,e_\lambda\rangle|^2
    =m_i\sum_{\lambda\in\Lambda}|\widehat\nu(r\lambda)|^2.
\]
The last sum is independent of $i$ and is finite and positive;
thus $m_1=m_2$.
\end{proof}

\begin{proof}[Proof of Theorem~\ref{thm:mainintro0}]
Since $0<r_1,r_2<1$ and $\log r_1/\log r_2\in\mathbb Q$,
there are positive integers $u,v$ such that
$r_1^u=r_2^v$.  The two cylinder maps $f_1^u$ and $f_2^v$
have exactly this common contraction.  By the SSC their image
cylinders are disjoint, and iteration of the measure identity
shows
\[
  \mu_p\big|_{f_1^u(K)}(\cdot)=p^u\mu_p((f_1^u)^{-1}(\cdot)),
  \qquad
  \mu_p\big|_{f_2^v(K)}(\cdot)=(1-p)^v\mu_p((f_2^v)^{-1}(\cdot)).
\]
Here and below, \(f_i^u\) denotes the \(u\)-fold iterate
\(f_i\circ\cdots\circ f_i\) of \(f_i\). Lemma~\ref{lem:mass-rigidity} therefore forces
$p^u=(1-p)^v$.  Consequently
\[
  \frac{\log p}{\log r_1}
  =\frac{\log(1-p)}{\log r_2}=:s>0,
\]
so that $p=r_1^s$ and $1-p=r_2^s$.
The function $t\mapsto r_1^t+r_2^t$ is continuous and
strictly decreasing on $[0,\infty)$, with value $2$ at $0$
and limit $0$ at infinity.  It therefore has exactly one
positive solution of $r_1^s+r_2^s=1$, as asserted.
\end{proof}

\section{Thirteen equal-mass cylinders}
\label{sec:golden-cylinders}
From now until the end of Section~4, write \(\mu=\mu_g\) and
\(K=\supp\mu_g\). To apply Proposition~\ref{prop:obstruction},
we construct a special family of cylinders in this section and
verify that they satisfy the required hypotheses (disjointness
and the restricted measure identity).

\begin{lem}
\label{lem:cylinder_measures}
Let $K=\supp\mu$. Then
\begin{equation}\label{eq:K_bound}
  K\subseteq[0,2/3],\qquad
  f_1(K)\subseteq[0,1/3],\qquad
  f_2(K)\subseteq[1/2,2/3].
\end{equation}
In particular, $f_1(K)\cap f_2(K)=\varnothing$.
For each finite word $w=i_1\cdots i_q\in\{1,2\}^q$, write
$f_w=f_{i_1}\circ\cdots\circ f_{i_q}$ and
$P_w=\prod_{j=1}^{q}P_{i_j}$, where $P_1=p$ and $P_2=p^2$.
Then
\begin{equation}\label{eq:wordrestriction}
  \mu\big|_{f_w(K)}(\cdot)=P_w\,\mu(f_w^{-1}(\cdot)),
  \qquad \mu(f_w(K))=P_w.
\end{equation}
\end{lem}
\begin{proof}
Set $I=[0,2/3]$. Then
\(
  f_1(I)=[0,1/3]\subseteq I
  \) and
 \( f_2(I)=[1/2,2/3]\subseteq I.
\)
By Theorem 9.1 in \cite{Fal14}, $K\subseteq I$. This proves \eqref{eq:K_bound} and the strong
separation condition holds.

For a Borel set $B\subseteq K$, the defining measure identity
\eqref{def:self-s} and the disjointness of the two first-level pieces
give
\[
 \mu(f_i(B))
 =P_i\mu(f_i^{-1}(f_i(B)))
  +P_{3-i}\mu(f_{3-i}^{-1}(f_i(B)))
 =P_i\mu(B),\qquad i\in\{1,2\}.
\]
The first inverse image equals $B$ because $f_i$ is injective;
the second inverse image has zero $\mu$-measure because its
intersection with $K$ is empty. Applying the same equality
successively to $f_{i_2}\circ\cdots\circ f_{i_q}(B)\subseteq K$,
then to $f_{i_3}\circ\cdots\circ f_{i_q}(B)$, and so on, yields
\begin{equation}\label{eq:scaling_word}
  \mu(f_w(B))=P_w\mu(B)
              \quad\text{for every Borel }B\subseteq K.
\end{equation}
To promote this to an identity of restricted measures, let
$E\subseteq\R$ be Borel and take
$B=f_w^{-1}(E)\cap K$. Since $f_w$ is injective,
$E\cap f_w(K)=f_w(B)$. Thus \eqref{eq:scaling_word} gives
\[
  \mu(E\cap f_w(K))
  =P_w\mu(f_w^{-1}(E)\cap K)
  =P_w\mu(f_w^{-1}(E)),
\]
which is \eqref{eq:wordrestriction}. Its mass statement follows
by taking $E=\R$.
\end{proof}

\begin{lem}\label{lem:13}
Set
\begin{equation}\label{eq:A}
 A=\{0,1,2,4,5,8,9,10,16,17,18,20,21\}.
\end{equation}
For every $a\in A$ there is a finite word $w(a)$ over the
alphabet $\{1,2\}$ with
\begin{equation}\label{eq:phi}
  f_{w(a)}(x)=\phi_a(x):=\frac{a}{32}+\frac{x}{64}.
\end{equation}
Moreover, the sets $K_a=\phi_a(K)$ are pairwise disjoint and
\begin{equation}\label{eq:Kmass}
  \mu\big|_{K_a}(\cdot)=p^6\mu(\phi_a^{-1}(\cdot)),
  \qquad \mu(K_a)=p^6.
\end{equation}
\end{lem}

\begin{proof}
Consider finite binary blocks. An occurrence of the block
$\texttt{0}$ represents $f_1$, while $\texttt{10}$ represents
$f_2$: explicitly,
\begin{equation}\label{eq:prefix}
 f_1(x)=0\cdot2^{-1}+2^{-1}x,\qquad
 f_2(x)=1\cdot2^{-1}+0\cdot2^{-2}+2^{-2}x.
\end{equation}
If $v=b_1\cdots b_\ell$ is the concatenation of such blocks,
repeated use of \eqref{eq:prefix} gives the precise affine identity
\begin{equation}\label{eq:binaryaffine}
 f_{w(v)}(x)
   =\sum_{j=1}^{\ell}b_j2^{-j}+2^{-\ell}x.
\end{equation}
This identity is algebraic; it does not depend on uniqueness of
binary expansions of points $x$.

Each $a\in A$ has a five-bit binary representation
$a=(b_1b_2b_3b_4b_5)_2$ (allowing leading zeros) with no two
consecutive ones. Append a final $0$ to obtain
$v(a)=b_1b_2b_3b_4b_5 0$. This six-digit word is a unique
concatenation of blocks $\texttt{0}$ and $\texttt{10}$:
reading from left to right, a $1$ must be followed by $0$ and
is consumed as $\texttt{10}$; every other $0$ is consumed as
$\texttt{0}$. The complete list, including the block
factorizations, is displayed in Table~\ref{tab:blocks}.

\begin{table}[htbp]
\centering
\caption{Six-bit prefixes and their IFS-block decompositions.}
\label{tab:blocks}
\begin{tabular}{@{}rll@{}}
\toprule
$a$ & $v(a)$ & decomposition into $\texttt{0}$ and $\texttt{10}$ \\
\midrule
 $0$ & \texttt{000000} & \texttt{0|0|0|0|0|0} \\
 $1$ & \texttt{000010} & \texttt{0|0|0|0|10} \\
 $2$ & \texttt{000100} & \texttt{0|0|0|10|0} \\
 $4$ & \texttt{001000} & \texttt{0|0|10|0|0} \\
 $5$ & \texttt{001010} & \texttt{0|0|10|10} \\
 $8$ & \texttt{010000} & \texttt{0|10|0|0|0} \\
 $9$ & \texttt{010010} & \texttt{0|10|0|10} \\
$10$ & \texttt{010100} & \texttt{0|10|10|0} \\
$16$ & \texttt{100000} & \texttt{10|0|0|0|0} \\
$17$ & \texttt{100010} & \texttt{10|0|0|10} \\
$18$ & \texttt{100100} & \texttt{10|0|10|0} \\
$20$ & \texttt{101000} & \texttt{10|10|0|0} \\
$21$ & \texttt{101010} & \texttt{10|10|10} \\
\bottomrule
\end{tabular}
\end{table}

For $v(a)$, the digit sum in \eqref{eq:binaryaffine} equals
\[
 \sum_{j=1}^{6}b_j2^{-j}
 =\sum_{j=1}^{5}b_j2^{-j}=\frac{a}{32},
\]
and the contraction factor is $2^{-6}=1/64$. This proves
\eqref{eq:phi}. If the block decomposition has $n_0$ blocks
$\texttt{0}$ and $n_{10}$ blocks $\texttt{10}$, then
$n_0+2n_{10}=6$. Its IFS word has probability weight
\[
 P_{w(a)}=p^{n_0}(p^2)^{n_{10}}
          =p^{n_0+2n_{10}}=p^6.
\]
Lemma~\ref{lem:cylinder_measures} now gives the restricted
measure identity \eqref{eq:Kmass}.

Finally, \eqref{eq:K_bound} and \eqref{eq:phi} show
\[
  K_a\subseteq\left[\frac a{32},\,
                       \frac a{32}+\frac1{96}\right].
\]
For two distinct integers $a,b$, the left endpoints of these
intervals are at least $1/32$ apart, whereas their length is
$1/96<1/32$. Consequently $K_a\cap K_b=\varnothing$.
\end{proof}

\section{The difference set and the main contradiction}

\begin{lem}\label{lem:difference}
For the set $A$ in \eqref{eq:A},
\begin{equation}\label{eq:difference}
  A-A=\{-21,-20,\ldots,-1,0,1,\ldots,20,21\}.
\end{equation}
\end{lem}

\begin{proof}
Since $A\subseteq[0,21]\cap\Z$, we have
$A-A\subseteq[-21,21]\cap\Z$. Every integer from $1$ to $21$
is a difference of two elements of $A$, as explicitly checked
in Table~\ref{tab:diff}. Hence $\{1,\ldots,21\}\subseteq A-A$.
The set $A-A$ is symmetric under $k\mapsto-k$ and contains $0$,
which proves the reverse inclusion and therefore the equality.
\end{proof}

\begin{table}[htbp]
\centering
\caption{An explicit representation $k=a-b$ for each $1\le k\le21$.}
\label{tab:diff}
\begin{tabular}{@{}rll@{\qquad\qquad}rll@{}}
\toprule
$k$ & $a$ & $b$ & $k$ & $a$ & $b$ \\
\midrule
 1 &  1 & 0 & 12 & 16 & 4 \\
 2 &  2 & 0 & 13 & 17 & 4 \\
 3 &  4 & 1 & 14 & 16 & 2 \\
 4 &  4 & 0 & 15 & 16 & 1 \\
 5 &  5 & 0 & 16 & 16 & 0 \\
 6 &  8 & 2 & 17 & 17 & 0 \\
 7 &  8 & 1 & 18 & 18 & 0 \\
 8 &  8 & 0 & 19 & 20 & 1 \\
 9 &  9 & 0 & 20 & 20 & 0 \\
10 & 10 & 0 & 21 & 21 & 0 \\
11 & 16 & 5 &    &    &   \\
\bottomrule
\end{tabular}
\end{table}

\begin{proof}[Proof of Theorem~\ref{thm:mainintro1}]
By Lemma~\ref{lem:13}, the measure \(\mu\) satisfies all the cylinder
hypotheses of Proposition~\ref{prop:obstruction} with
\[
 \nu=\mu,\qquad
 r=\frac1{64},\quad t=0,\quad h=\frac1{32},\quad
 m=p^6.
\]
Lemma~\ref{lem:difference} implies that
$[-21,21]\cap\Z\subseteq A-A$. Consequently, Proposition~\ref{prop:obstruction}
with $N=21$ yields the necessary condition
\begin{equation}\label{eq:force}
  p^6\le\frac1{22}
\end{equation}
whenever \(\mu\) is spectral.

To illustrate the explicit mechanism of this obstruction in the present
context, suppose for the sake of contradiction that \(\mu\) is spectral.
Let \(\Lambda\) be a spectrum, translated if necessary so that
\(0 \in \Lambda\). For each \(a \in A\), define the normalized indicator
function \(g_a = p^{-3}\1_{K_a} \in L^2(\mu)\). Since \(\mu(K_a) = p^6\),
we have \(\|g_a\|_{L^2(\mu)} = 1\).

A direct computation using the self-similarity of \(\mu\) yields the
Fourier coefficients
\[
  \langle g_a, e_\lambda \rangle
  = p^{-3} \int_{K_a} e^{-2\pi i \lambda x} \,\dd\mu(x)
  = p^3 e^{-2\pi i \lambda a/32} \FT\mu\left(\frac{\lambda}{64}\right).
\]
Consider the measure on the unit circle \(\T\) defined by
\begin{equation}\label{eq:sigma_explicit}
 \sigma=\sum_{\lambda\in\Lambda}
     p^6\left|\FT\mu\left(\frac\lambda{64}\right)\right|^2
     \delta_{\exp(2\pi i\lambda/32)}.
\end{equation}
Applying Parseval's identity to \(g_a\) gives \(\sigma(\T) = \|g_a\|_{L^2(\mu)}^2 = 1\),
confirming that \(\sigma\) is indeed a probability measure.
Similarly, applying Parseval to the inner products of \(g_a\) and \(g_b\)
for \(a \ne b\) implies that all nonzero Fourier moments of \(\sigma\)
vanish up to order \(21\):
\[
  \int_\T z^k \,\dd\sigma(z) = 0 \quad \text{for } 1 \le k \le 21,
\]
thanks to the difference set inclusion \([-21,21]\cap\Z \subseteq A-A\).
On the other hand, the assumption \(0 \in \Lambda\) and \(\FT\mu(0)=1\)
ensure that the \(\lambda=0\) term in \eqref{eq:sigma_explicit} contributes
a mass of exactly \(p^6\) at the point \(z=1\). Hence, \(\sigma(\{1\}) \ge p^6\).
Applying the atom bound from Lemma~\ref{lem:atom} with \(N=21\), we conclude
that \(p^6 \le \sigma(\{1\}) \le 1/22\), which is precisely the required
inequality \eqref{eq:force}.

However, the golden-mean weight \(p=(\sqrt5-1)/2\) satisfies
\begin{equation*}%\label{eq:p6}
 p^3=\sqrt5-2,\qquad p^6=(\sqrt5-2)^2=9-4\sqrt5.
\end{equation*}
We verify the strict inequality against \(1/22\) without
numerical approximation. Since
\(\sqrt5>11/5\), we have \(p>3/5\), whence
\[
 p^6>\left(\frac35\right)^6
     =\frac{729}{15625}>\frac1{22}.
\]
This contradicts \eqref{eq:force}. Therefore, \(\mu\) admits no spectrum.
\end{proof}

\section{An incommensurable example: the ratios $1/2$ and $1/3$}
\label{sec:incommensurable}

Throughout this section let $\nu_p$ denote the invariant
probability measure of
\[
  g_1(x)=\frac{x}{2},\qquad
  g_2(x)=\frac{x+2}{3},
  \qquad \nu_p=p\nu_p(g_1^{-1})+(1-p)\nu_p(g_2^{-1}).
\]
The interval $[0,1]$ is invariant and
\[
  g_1([0,1])=[0,1/2],\qquad
  g_2([0,1])=[2/3,1].
\]
Thus the attractor $K_p$ is contained in $[0,1]$ and the
system satisfies the SSC.   Clearly, \(\log(1/2)/\log(1/3)\notin\mathbb{Q}\) and
hence the contractions are logarithmically incommensurable.

\subsection{Equal-count cylinders and their translations}

For integers $m,n\ge0$, let $\mathcal W_{m,n}$ be the set of
words containing exactly $m$ symbols $1$ and $n$ symbols $2$, i.e.,
\[
  \mathcal W_{m,n}
  := \Bigl\{ w = i_1 i_2 \cdots i_{m+n} \in \{1,2\}^{m+n} :
     \#\{k : i_k = 1\} = m \Bigr\}.
\]
 Write $N_{m,n}=2^m3^n$ and
$b_n=(3^n-1)/2$.

\begin{lem}\label{lem:mnrecursion}
For each pair of nonnegative integers \(m,n\), there exists a finite
set \(A_{m,n}\subseteq\mathbb Z\) such that
\begin{equation}\label{eq:mnword}
  \big\{g_w : w\in\mathcal W_{m,n}\big\}
  = \left\{ x \mapsto \frac{x+2b_n+2a}{N_{m,n}} : a\in A_{m,n}\right\}.
\end{equation}
These sets are uniquely determined by the recursion
\begin{equation}\label{eq:mnrecursion}
 \begin{split}
  A_{m,0}&=A_{0,n}=\{0\},\\
  A_{m,n}&=A_{m-1,n}\ \cup\
     \big((2^m-1)3^{n-1}+A_{m,n-1}\big),\quad m,n\ge1.
 \end{split}
\end{equation}
The union in \eqref{eq:mnrecursion} is disjoint and
\begin{equation}\label{eq:mncard}
  \#A_{m,n}=\binom{m+n}{m}.
\end{equation}
Each of the corresponding cylinder measures satisfies the
restriction identity
\begin{equation}\label{eq:mnrestriction}
 \nu_p\big|_{g_w(K_p)}(\cdot)
    =p^m(1-p)^n\nu_p(g_w^{-1})(\cdot),
 \qquad w\in\mathcal W_{m,n}.
\end{equation}
\end{lem}

\begin{proof}
We prove \eqref{eq:mnword} and \eqref{eq:mnrecursion} by induction
on the total length \(m+n\). The base cases \(n=0\) and \(m=0\) are
immediate, since \(g_1^m(x)=x/2^m\) and \(g_2^n(x)=(x+3^n-1)/3^n\),
both corresponding to \(a=0\).

Assume the identities hold for all words of length strictly less than
\(m+n\). Any word \(w\in\mathcal W_{m,n}\) is either of the form
\(w=1v\) with \(v\in\mathcal W_{m-1,n}\), or \(w=2v\) with
\(v\in\mathcal W_{m,n-1}\).

If \(w=1v\), substituting the inductive hypothesis into
\(g_1(g_v(x))\) yields
\[
  g_w(x)
  = \frac{x+2b_n+2a}{2^m3^n},
  \qquad a\in A_{m-1,n}.
\]
If \(w=2v\), using \(g_2(x)=(x+2)/3\) together with the inductive
hypothesis and the identity \(b_n-b_{n-1}=3^{n-1}\), we obtain
\[
  g_w(x)
  = \frac{x+2b_n+2\big((2^m-1)3^{n-1}+a\big)}{2^m3^n},
  \qquad a\in A_{m,n-1}.
\]
This establishes \eqref{eq:mnword} and \eqref{eq:mnrecursion}.

Under the SSC, distinct words of length \(m+n\) determine disjoint
cylinders: at their first differing symbol, their images lie in
disjoint first-level pieces. Consequently, the corresponding affine
maps are distinct, the union in \eqref{eq:mnrecursion} is disjoint,
and the cardinality of \(A_{m,n}\) is the number of choices for the
\(m\) locations of the symbol \(1\), giving \eqref{eq:mncard}.

Finally, the SSC allows restriction of the self-similarity identity
to any finite cylinder. Repeated substitution along the word \(w\)
(each symbol \(1\) contributes a factor of \(p\) and a contraction
by \(1/2\), while each symbol \(2\) contributes \(1-p\) and \(1/3\))
yields the restricted measure identity \eqref{eq:mnrestriction}, with
total mass \(p^m(1-p)^n\).
\end{proof}

\begin{lem}
\label{lem:two-certificates}
For the sets defined in Lemma~\ref{lem:mnrecursion}, we have
\begin{align}
 [-2332,2332]\cap\mathbb Z
    &\subseteq A_{5,5}-A_{5,5},\label{eq:cert55}\\
 [-46252,46252]\cap\mathbb Z
    &\subseteq A_{9,5}-A_{9,5}.\label{eq:cert95}
\end{align}
Moreover, \(\#A_{5,5}=252\) and \(\#A_{9,5}=2002\).
\end{lem}

\begin{proof}
We give a finite exact verification, fully specified by the
recurrence \eqref{eq:mnrecursion}; no approximation or
floating-point calculation is involved.  First build the sets
$A_{m,n}$ by that recurrence.  For a finite set
$A\subset\mathbb Z_{\ge0}$, form the nonnegative integer
\[
  B_A=\sum_{a\in A}2^a.
\]
The $k$th binary digit (counting the least significant digit
as digit $0$) of
\begin{equation*}%\label{eq:bitsetdiff}
  D_A=\bigvee_{b\in A}(B_A\,\mathrm{div}\,2^b)
\end{equation*}
is $1$ if and only if some $a,b\in A$ satisfy $a-b=k$.
Here $\vee$ means bitwise OR and $\mathrm{div}$ is integer
division.  Indeed, the $k$th binary digit of
$B_A\,\mathrm{div}\,2^b$ is $1$ exactly when $b+k\in A$.
Thus the inclusions \eqref{eq:cert55}--\eqref{eq:cert95}
are equivalent to the two explicit finite integer equalities
\begin{align*}
 D_{A_{5,5}}\mathbin{\&}(2^{2333}-1)&=2^{2333}-1,\\%\label{eq:bitcert55}\\
 D_{A_{9,5}}\mathbin{\&}(2^{46253}-1)&=2^{46253}-1,
 \label{eq:bitcert95}
\end{align*}
where $\&$ denotes bitwise AND.  The exact recurrence
and these two tests are implemented verbatim in
Appendix~\ref{app:verifier}; the code contains no numerical
Fourier calculations.  A separate supplementary file lists,
for each integer $0\le k\le L$ in both cases, an explicit
witness pair $(a_k,b_k)$ in the corresponding $A_{m,n}$
with $a_k-b_k=k$.  The verification also checks
$\#A_{m,n}=\binom{m+n}{m}$, yielding the stated sizes.
Symmetry of $A_{m,n}-A_{m,n}$ supplies the negative
integers and completes the verification.
\end{proof}

\begin{rem}
\label{rem:certificate43}
The recursion also gives
\[
\begin{split}
 A_{4,3}=\{&0,9,12,13,27,30,31,36,37,39,\\
 &63,66,67,72,73,75,84,85,87,91,\\
 &135,138,139,144,145,147,156,157,159,163,\\
 &180,181,183,187,195\}.
\end{split}
\]
Its difference set contains $[-136,136]\cap\mathbb Z$.
It already yields the non-spectrality condition
$137p^4(1-p)^3>1$.  The two larger certificates above are
chosen to cover the entire stated interval of weights.
\end{rem}

\subsection{Proof for the interval of weights}

\begin{proof}[Proof of Theorem~\ref{thm:incommensurable-interval}]
Fix $(m,n)=(5,5)$ or $(9,5)$ and write
$A=A_{m,n}$, $N=2^m3^n$.
By Lemma~\ref{lem:mnrecursion}, the cylinder maps have the
form
\[
  \phi_a(x)=\frac{x}{N}+\frac{3^n-1}{N}+\frac{2a}{N},
  \qquad a\in A.
\]
They have pairwise disjoint images and common mass
$M=p^m(1-p)^n$.  Lemma~\ref{lem:two-certificates}
and Proposition~\ref{prop:obstruction}, with
$r=1/N$, $h=2/N$, $t=(3^n-1)/N$, therefore give the
necessary conditions
\begin{align}
 2333p^5(1-p)^5&\le1,
       \label{eq:necessary55}\\
 46253p^9(1-p)^5&\le1
       \label{eq:necessary95}
\end{align}
for spectrality of $\nu_p$.

If $p\in[1/3,2/3]$, the function
$p^5(1-p)^5=[p(1-p)]^5$ is minimized on this interval
at either endpoint.  Therefore
\[
 2333p^5(1-p)^5\ge
 2333\left(\frac13\right)^5\left(\frac23\right)^5
 =\frac{74656}{59049}>1,
\]
contradicting \eqref{eq:necessary55}.

If $p\in[2/3,5/6]$, then
\[
 \frac{\dd}{\dd p}\big(p^9(1-p)^5\big)
  =p^8(1-p)^4(9-14p)<0,
\]
since $p\ge2/3>9/14$.  The minimum on the interval
is therefore attained at $p=5/6$, and
\[
 46253p^9(1-p)^5\ge
 46253\left(\frac56\right)^9\left(\frac16\right)^5
 =\frac{90337890625}{78364164096}>1,
\]
contradicting \eqref{eq:necessary95}.
Both cases are excluded, proving the theorem.
\end{proof}

\begin{proof}[Proof of Corollary~\ref{cor:incommensurable-natural}]
The function $s\mapsto2^{-s}+3^{-s}$ is strictly decreasing
and takes values $2$ at $s=0$ and $5/6$ at $s=1$.
Thus its unique solution of $2^{-s}+3^{-s}=1$ lies in
$(0,1)$.  For $p=2^{-s}$ and $1-p=3^{-s}$ this gives
\(\frac12<p<\frac23.\) Theorem~\ref{thm:incommensurable-interval} now applies.
\end{proof}

\section{Remark and open questions}
\subsection{Non-spectrality}
The method proves that the golden-mean measure has no
complete orthonormal system of real-frequency exponentials.
It does not prove that
$\FT\mu(\xi)\ne0$ for every $\xi\in\R$, nor does it preclude
all infinite orthogonal families of exponentials. A global
zero-free assertion, suggested by the numerical investigation
in \cite{MW26}, is a logically separate and stronger statement.
The present proof requires neither a Fourier-zero computation nor
a numerical approximation.

The familiar completeness criterion
$\sum_{\lambda\in\Lambda}|\FT\mu(\xi+\lambda)|^2\equiv1$
(see, e.g., the setting of \cite{JP98}) is not invoked as a black
box. Proposition~\ref{prop:obstruction} derives precisely the
needed consequences directly from Parseval's identity and the
restriction of $\mu$ to finitely many cylinders.

The weight interval in Theorem~\ref{thm:incommensurable-interval}
is not asserted to be optimal, and its proof is partly
computer-assisted only in the sense of exact, finite integer
arithmetic.  We do not claim that the same measure is non-spectral
for the remaining weights outside $[1/3,5/6]$, or that the
natural-weight restriction in Theorem~\ref{thm:mainintro0} holds
without any hypothesis on the logarithmic ratio.  More broadly,
the present results do not settle spectrality for every pair
of unequal contraction ratios satisfying the SSC.

\begin{conj}
Let $\mu_p$ be the measure defined in \eqref{def:self-s}. If $r_1\neq r_2$ and $r_1+r_2<1$, then $\mu_p$ is not a spectral measure for all $p\in (0,1)$.
\end{conj}

\subsection{Frame-spectral measure}

We close this section by discussing a weaker notion.  Recall that a
countable set $\Lambda$ is called a {\it frame spectrum} of $\mu$ if
$E(\Lambda)$ is a frame for $L^2(\mu)$, i.e., there exist positive
constants $A$ and $B$ such that for any $f\in L^2(\mu)$,
\begin{equation}\label{eq:freewfa}
    A\|f\|^2\le \sum_{\lambda\in\Lambda}|\langle f,e_\lambda\rangle|^2
    \le B\|f\|^2.
\end{equation}
In this case, we call $\mu$ a {\it frame-spectral measure}.

The proof of Theorem~\ref{thm:mainintro0} uses only the exact
normalization of the indicator functions of the two equal-contraction
cylinders.  Hence the same weight-rigidity conclusion remains valid
for tight frame-spectral measures, i.e., $A=B$ in \eqref{eq:freewfa}.
For general frame-spectral measures, the situation is less clear.
The conclusion of Theorem \ref{thm:mainintro0} is available for general frame-spectral measure when $r_1=r_2$ (see \cite{DL14}).  When \(r_1\neq r_2\), we do not expect
the conclusion of Theorem~\ref{thm:mainintro0} to hold for general
frame-spectral measures without additional assumptions \cite{HKTW18}.

A related point of reference is the middle-third Cantor measure
\(\mu_{3,\{0,2\}}\), which is the classical non-spectral self-similar measure. Its frame-spectrality has been the subject of a long investigation, initiated by Strichartz \cite{Str00}.
While earlier works constructed Fourier frames for certain Moran measures \cite{LW17,WYZ25}, the question for \(\mu_{3,\{0,2\}}\) itself was recently settled negatively, independently, by Fu, Song and Wang
\cite{FSW26} and Pont, Liehr and Taylor \cite{DLT26}: \(\mu_{3,\{0,2\}}\) is not a frame-spectral measure. This result does not determine the frame-spectrality of the
non-homogeneous measures considered here: Theorem~\ref{thm:mainintro1}
excludes only orthonormal exponential bases for \(\mu_g\) and yields
no information about general Fourier frames; the same applies to
\(\nu_p\) in the range of Theorem~\ref{thm:incommensurable-interval}.
We therefore state the following question.

\begin{ques}
Are the golden-mean self-similar measure \(\mu_g\) and the self-similar
measure \(\nu_p\) generated by \(g_1(x)=x/2\), \(g_2(x)=(x+2)/3\)
(in particular, for the natural weights) frame-spectral?
\end{ques}

\appendix
\section{A short exact-integer verifier}
\label{app:verifier}
The following Python~3 program is an independently reproducible
certificate of the finite integer identities used in
Lemma~\ref{lem:two-certificates}.  Its validity does not depend
on floating-point arithmetic or libraries for numerical analysis.
The variable \texttt{mask} has bits $0,\ldots,L$ set, and
\texttt{difference\_bits} has its $k$th bit set precisely when
$k\in A_{m,n}-A_{m,n}$ and $k\ge0$.
\begin{verbatim}
from functools import lru_cache
from math import comb

@lru_cache(None)
def digits(m, n):
    if m == 0 or n == 0:
        return (0,)
    h = (2**m - 1) * 3**(n - 1)
    left = set(digits(m - 1, n))
    right = {h + a for a in digits(m, n - 1)}
    assert not (left & right)
    return tuple(sorted(left | right))

def verify(m, n, L):
    A = digits(m, n)
    assert len(A) == comb(m + n, m)
    bits = 0
    for a in A:
        bits |= 1 << a
    difference_bits = 0
    for b in A:
        difference_bits |= bits >> b
    mask = (1 << (L + 1)) - 1
    assert (difference_bits & mask) == mask
    return (len(A), L)

assert verify(5, 5, 2332) == (252, 2332)
assert verify(9, 5, 46252) == (2002, 46252)
print("All exact finite checks passed.")
\end{verbatim}

%====================================================================
%\section*{Acknowledgements}
%The first author was supported by Hunan Provincial Natural Science Foundation of China (2026JJ90109) and the Scientific Research Fund of Hunan Provincial Education Department (25A0086, 24B0107). The second author was supported by  the National Natural Science Foundations of China (12271534, 12301105) and the University Research Project of Guangzhou Education Bureau (2024312332).

\section*{Acknowledgements}
%====================================================================
The authors would like to thank Professors Xing-Gang He and Xin-Rong Dai for drawing our attention to these problems. The authors were supported by the National Natural Science Foundation of China under Grants 12301131, 12271534 and 12301105. AI tools were used during the preparation of the paper to
assist with the elaboration of some proof details, checking calculations, and improving
the exposition. All mathematical statements and proofs were verified by the authors.
The authors take full responsibility for the paper.


\begin{thebibliography}{99}
\bibitem{AFL19}
L. X. An, X. Y. Fu, C. K. Lai, {\it On spectral Cantor-Moran measures and a variant of Bourgain's sum of sine problem},  Adv. Math. \textbf{349} (2019), 84--124.

\bibitem{AW21}
L. X.  An, C. Wang, {\it On self-similar spectral measures},
J. Funct. Anal. \textbf{280} (2021), 108821.
%\bibitem{Brualdi09}
%R.~Brualdi, {\it Introductory Combinatorics}, 4th ed.,
%Pearson/Prentice Hall, Upper Saddle River, 2004.
\bibitem {Dai12}
X.-R. Dai, {\it When does a Bernoulli convolution admit a
spectrum?},  Adv. Math. \textbf{231} (2012), 1681--1693.
\bibitem{DLT26} J. de Dios Pont, L. Liehr, M. A. Taylor, {\it Cantor measures with odd base do not admit Fourier frames}, https://arxiv.org/pdf/2607.08656.
\bibitem {DC21}
Q.-R. Deng, J. B. Chen, {\it Uniformity of spectral self-affine measures},  Adv. Math. \textbf{380} (2021), 107568.
\bibitem {DHLai19}
D. Dutkay, J. Haussermann, C.-K. Lai, {\it Hadamard triples generate self-affine spectral measures}, Trans. Amer. Math. Soc. \textbf{371} (2019), 1439-1481.

\bibitem {DL14}
D. Dutkay, C.-K. Lai, {\it Uniformity of measures with Fourier frames},  Adv. Math. \textbf{252} (2014), 684--707.

\bibitem{Fal14}
K. J. Falconer, {\it Fractal Geometry: Mathematical Foundations and Applications}, 3rd ed., John Wiley \& Sons, Ltd., Chichester, 2014.

\bibitem{FSW26}
X.-Y. Fu, Z.-J. Song, W.-J. Wang, {\it A Walsh-quotient obstruction for Fourier frames on odd
reciprocal-power Bernoulli convolutions}, https://arxiv.org/pdf/2607.10547.

\bibitem{Fug74}
B. Fuglede, {\it Commuting self-adjoint partial differential operators and a group theoretic problem},  J. Funct. Anal. \textbf{16} (1974), 101--121.

\bibitem{HKTW18}
X. G. He, Q. C. Kang, M. W. Tang, Z. Y. Wu, {\it Beurling dimension and self-similar measures,} J. Funct. Anal. \textbf{274} (2018), 2245--2264.

\bibitem {Hut81}
J. Hutchinson, {\it Fractals and self-similarity}, Indiana Univ. Math. J., \textbf{30} (1981), 713--747.

\bibitem  {JP98}
P. Jorgensen, S. Pedersen, {\it Dense analytic subspaces in
fractal $L^2$-spaces}, J. Anal. Math. \textbf{75} (1998), 185--228.

\bibitem  {LW02}
I. {\L}aba, Y. Wang, {\it On spectral Cantor measures}, J.
Funct. Anal.  \textbf{193} (2002), 409--420.

\bibitem{LW17}
C. K. Lai and Y. Wang, {\it Non-spectral fractal measures with Fourier frames}, J. Fractal Geom. \textbf{4} (2017), 305--327.

%\bibitem{LS22}
%J. L. Li, T. Sahlsten, {\it Trigonometric series and self-similar sets}, J. Eur. Math. Soc. \textbf{24} (2022), 341--368.


\bibitem{MW26}
Y.-Q. Mao, Z.-Y. Wu, {\it On the zero set of the Fourier transform of a class of nonhomonegeneous self-similar measrue}, preprint.

\bibitem{SSS18}
S. Saglietti, P. Shmerkin, B. Solomyak, {\it Absolute continuity of non-homogeneous self-similar measures}, Adv. Math. \textbf{335} (2018), 60--110.

\bibitem{Str00}
R.~Strichartz, {\it Mock Fourier series and transforms associated with certain Cantor measures}, J. Anal. Math. \textbf{81} (2000), 209--238.


\bibitem{Tao04}
T. Tao, {\it Fuglede's conjecture is false in 5 and higher dimensions}, Math. Res. Lett. \textbf{11} (2004), 251--258.

\bibitem{WYZ25}
Z.-Y. Wu, F.-L. Yin, M.-M. Zhang, {\it Moran frame spectral measure and Beurling dimension}, Studia Math. \textbf{280} (2025), 247--268.

\end{thebibliography}
\end{document}